\documentclass[preprint,12pt]{elsarticle} % 1 interval

\usepackage[T2A]{fontenc}
\usepackage[cp866]{inputenc}
\usepackage{graphicx}
\usepackage{epsfig}
\usepackage{amssymb, amsmath}
\usepackage{color}

\everymath{\displaystyle}

\newtheorem{theorem}{Theorem}

\journal{EJOR}

\begin{document}

\begin{frontmatter}



\title{Random sampling: Billiard Walk algorithm\tnoteref{label1}}
\tnotetext[label1]{The work was supported by Laboratory of
Structural Methods of Data Analysis in Predictive Modeling in Moscow
Institute of Physics and Technology (``mega-grant'' of the Russian
Government) and by RFFI grant 13-07-12111 ofi-m.}


\author[EG]{Elena Gryazina\corref{cor1}}
\ead{gryazina@gmail.com}
\cortext[cor1]{Corresponding author}
\address[EG]{Institute for Control Sciences RAS, Moscow, Russia; +7 495 334
8829}

\author[BP]{Boris Polyak}
 \address[BP]{Institute for Control Sciences RAS and Laboratory of Structural Methods
 of Data Analysis in Predictive Modeling in Moscow
Institute of Physics and Technology, Moscow, Russia}

\begin{abstract}
Hit-and-Run is known to be one of the best random sampling
algorithms, its mixing time is polynomial in dimension.
However in practice, the number of steps required to obtain
uniformly distributed samples is rather high. We propose a new random
walk algorithm based on billiard trajectories. Numerical experiments
demonstrate much faster convergence to the uniform distribution.
\end{abstract}


\end{frontmatter}

\section{Introduction}

%\cite{Sobol99}

% about applications of sampling
Generating points uniformly distributed in an arbitrary bounded
region $Q \subset \mathbb{R}^n$ finds applications in
many computational problems \cite{TeCaDa:04,Rubin_MCbook:2008}.

Straightforward sampling techniques are usually based on one of the
three approaches: rejection, transformation, and composition. In the
rejection approach, the region of interest $Q$ is embedded into a
region with available uniform sampler $B$ (usually a box or a
ball). At the next step, samples that do not belong to $Q$ are
rejected. Assume~$Q$ is the unit ball, and the bounding region $B$ is
the box $[-1, 1]^n$. Then for $n = 2k$, the ratio of the volumes of
the box and the ball is equal to
$q=\frac{\text{Vol}(Q)}{\text{Vol}(B)} = \frac{\pi ^k}{k! 2^k},$
thus  $ q\approx 10^{-8}$ for $n=20$, so that one has to generate
$\sim 10^8$ samples to obtain just a few of them in $Q$. For
polytopes this ratio can be much smaller. Another way to exploit
pseudo-random number generator for a simple region $B$ is to map $B$
onto $Q$ via a smooth deterministic function with constant Jacobian.
For instance, to obtain uniform samples in $Q=\{x: x^T A x \le 1\}$,
$A$ being a positive definite matrix, it suffices to generate samples
$y$ uniformly in the unit ball $||y||_2 \le 1$ and transform them as
$x = A ^{-1/2} y$. Unfortunately, such a transformation exists just
for a limited class of regions. In the composition approach, the
set~$Q$ is partitioned into a finite number of sets that can be
efficiently sampled. For instance, a polytope can be partitioned
into  simplices, but the large number of them makes the procedure
computationally hard.

Other sampling procedures use modern versions of the Monte Carlo
technique based on the Markov Chain Monte Carlo (MCMC) approach
\cite{Gilks_book96,Diaconis09}. For instance, efficient algorithms
for computing volumes using random walks can be found in
\cite{DyerFriKannan:91,LovaszSomonovits:93,LovaszDeak:12}. One of
the most famous and efficient algorithms of the MCMC type is
Hit-and-Run (HR), which was originally proposed by Turchin
\cite{Turchin_eng} and independently by Smith \cite{Smith}. The
brief description of the HR algorithm is as follows. At every step
HR generates a random direction uniformly over the unit sphere and
picks the next point uniformly on the segment of the straight line
in the given direction in $Q$. HR is applicable to various control
and optimization problems
\cite{PolGry:IFAC08,PolGry_AOR:09,DaShPo:SIAM10} as well as to
simulation-based multiple criteria decision analysis
 \cite{Tervonen_etal:EJOR:2013}. Unfortunately, even for simple  ``bad'' sets,
such as level sets of ill-posed functions, HR techniques fail or
become computationally inefficient.

A variety of applications and drawbacks of the existing techniques
provides much room for improving and developing new sampling
algorithms. For instance, there were attempts to exploit the
approach proposed for interior-point methods of convex optimization
\cite{NestNem_book94} and to combine it with MCMC algorithms. As a
result, the Barrier Monte Carlo method \cite{PolGry:MSC10} generates
random points with better uniformity properties as compared to the
standard Hit-and-Run. On the other hand, the complexity of every
iteration is in general high enough (the calculation of
$\left(\nabla^2 F(x)\right)^{-1/2}$ is required, where $F(x)$ is a
barrier function of the set). Moreover, the Barrier Monte Carlo
method does not accelerate convergence for simplex-like sets.

% about billiards (and new method)

In this paper we propose a new random walk algorithm motivated by
physical phenomena of gas diffusing in a vessel. A particle of gas
moves with a constant speed until it meets the boundary of the vessel,
then it reflects (the angle of incidence equals the angle of
reflection) and so on. When a particle hits another one, its
direction and speed change. In our simplified model we assume that
the direction changes randomly, while the speed remains the same. Thus our
model combines the ideas of the Hit-and-Run technique and use of the billiard
trajectories.  There exists a vast literature on mathematical
billiards, and many useful facts can be extracted from there
\cite{Tabachnikov:95, GalpZem:90_billiards, Sinai:70, Sinai:78,
Kozlov:billiards_book:91}. The traditional theory addresses the behavior
of one particular billiard trajectory in different billiard tables,
their ergodic properties, and the conditions for the existence of
periodic orbits. In stochastic analogs of the classical billiard
\cite{evans2001stochastic}, a direction after reflection is
chosen randomly uniformly. Shake-and-Bake algorithms are based on
stochastic billiards and generate points on the boundary of a convex set
\cite{Shake-n-Bake:91}. The recently proposed version of the Shake-and-Bake
algorithm \cite{dieker2013StochBill} exhibits polynomial-time
convergence to the uniform distribution. Our algorithm is aimed at sampling
the interior of a set (actually, later in the text we consider open regions).
Besides that, we extend billiard trajectories of random length
keeping the standard reflection law. Such an incorporation of randomness
also improves the ergodic properties.
% The paper is organized...

The paper is organized as follows. In Section 2 we present a novel
sampling algorithm and prove that it produces asymptotically
uniformly distributed samples in $Q$. In Section 3 we pay much
attention to some properties of the Billiard Walk (BW),
implementation issues are discussed as well. Simulation of BW for
particular test domains is presented in Section 4. Much attention is
devoted to the capability of BW to get out of the corner, in comparison
with HR. Here we consider just the most demonstrative types of
geometry. In Section 5 we briefly discuss possible applications of
the algorithm.

\section{Algorithm}

Assume there is a bounded, open connected set $Q\subset
\mathbb{R}^{n}$, $n\geq 2$, and a point $x^0\in Q$. Our aim is to generate
asymptotically uniform samples $x^i \in Q$, $i = 1, \dots, N$.

\begin{figure}
\begin{center}
\includegraphics[width=8.4cm]{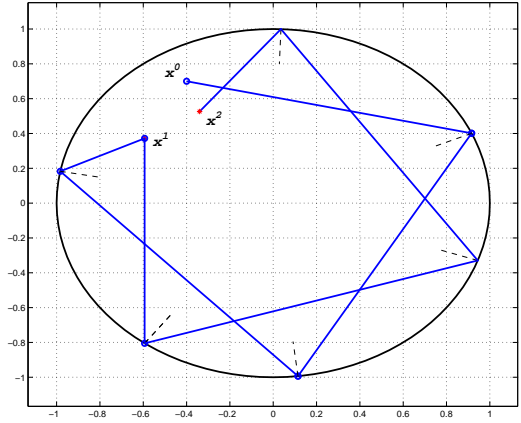}
\caption{Billiard Walk.}\label{fig:HRR}
\end{center}
\end{figure}

The new BW algorithm generates a random direction uniformly over the
unit sphere. Then the next sample is chosen as the end-point of the
billiard trajectory of length $\ell$. This length is chosen
randomly; i.e., we assume that the probability of collision with another
particle is proportional to $\delta t$ for small time instances
$\delta t$, this validates the formula for $\ell$ in the algorithm
below. The scheme of the method is given in Fig. \ref{fig:HRR}, while
the precise routine is as follows.

{\bf Algorithm of Billiard Walk (BW).}
\begin{enumerate}

\item[1.]
Take $x^0\in Q$; $i=0$, $x=x^0$.

\item[2.] Generate the length of the trajectory
$\ell = - \tau \text{log } \xi$, $\xi$ being uniform random on
$[0,1]$, $\tau$ is a specified constant parameter of the algorithm.

\item[3.]
Pick a random direction $d\in\mathbb{R}^{n}$ uniformly distributed
over the unit sphere (i.e., $d=\xi/\|\xi\|$, where
$\xi\in\mathbb{R}^{n}$ has the standard Gaussian distribution).
Construct a billiard trajectory starting at $x^i$ and having initial
direction $d$. When the trajectory meets the boundary with internal
normal $s$, $||s||=1$, the direction is changed as
$$
d \rightarrow d - 2(d,s)s,
$$
where $(d,s)$ is the scalar product.

\item[4.]
If a point with nonsmooth boundary is met or the number of
reflections exceeds $R$, go to step 3. Otherwise proceed until the
length of the trajectory equals $\ell$.

\item[5.] $i = i+1$, take the end-point as $x^{i+1}$ and go to step 2.
\end{enumerate}

The algorithm involves two parameters $\tau$ and $R$ and we discuss
their choice below.

We prove asymptotical uniformity of the samples produced by BW for
convex and nonconvex cases separately. The requirements on $Q$ are
different for these two cases, while the sampling algorithm remains
the same. Consider the Markov Chain induced by the BW algorithm
$x^0$, $x^1, \dots$. For an arbitrary measurable set $A\subseteq Q$,
denote by $\mathbf{P}(A|x)$ the probability of obtaining $x^{i+1}\in
A$ for $x^i = x$ by the BW algorithm. Then $\mathbf{P}_N(A|x)$ is
the probability to get $x^{i+N} \in A$ for $x^i = x$. We also denote
by $p(y|x)$ the probability density function for $\mathbf{P}(A|x)$,
i.e. $\mathbf{P}(A|x) = \int\limits_A p(y|x) dy$.

\begin{theorem}
Assume $Q$ is an open bounded convex set in $\mathbb{R}^n$,
the boundary of  $Q$ is piecewise smooth. Then the
distribution of points $x^i$ generated by the BW algorithm tends to the
uniform one over $Q$, i.e.
$$\lim_{N\rightarrow \infty} \mathbf{P}_N(A|x)  = \lambda (A) $$
for any measurable $A\subseteq Q$, $\lambda(A) = \text{Vol}(A) / \text{Vol}(Q)$
and any starting point $x$.
\end{theorem}

\textit{Proof.} First, the algorithm is well-defined: at step 4 with
zero probability the algorithm sticks at a point with nonsmooth
boundary. On the other hand  $\ell$ and $d$ are chosen in such a way that,
with positive probability, $x^{i+1}$ is obtained by less than $R$
reflections (see detailed discussion of ``bad'' situations in
Subsections~3.1 and~3.2).

In view of Theorem 2 in \cite{Smith} based on the asymptotic
properties of Markov Chains, the two assumptions on $p(y|x)$ imply
that the uniform distribution over $Q$ is a unique stationary
distribution, and it is achieved for any starting point $x\in Q$.
The first assumption requires the existence of $p(y|x)$ and its
symmetry; the second assumption claims its positivity $p(y|x)>0$ for
all $x,y \in Q$.

Now we show that there exists a probability density function; i.e. for any
$x,y \in Q$, the transition probability from $x$ to a small
neighborhood $\delta y$ of $y$ is proportional to the volume of
$\delta y$. Among the trajectories proceeding from $x$ to $\delta
y$, there exist a conic bundle of trajectories with no reflections,
as well as some trajectories with $1, 2, \dots, R$ reflections. For
a bundle of trajectories with no reflections $\mathbf{P}(\delta
y|x)\sim \mathbf{P}(\delta \theta)\mathbf{P}(\delta \ell)$, where
$\mathbf{P}(\delta \theta)$ is the probability of choosing the
spatial angle and $\mathbf{P}(\delta \ell)$ is the probability of
choosing a certain trajectory length $\ell\in \delta\ell$ while
$p\sim q$ means ``$p$ is proportional to $q$'' . $\mathbf{P}(\delta
\theta)$ is proportional to the volume of the base of the cone,
$\mathbf{P}(\delta \ell) \sim \delta \ell$, thus $\mathbf{P}(\delta
y|x) \sim \text{vol}(\delta y)$.

The bundles of trajectories with reflections are also cones with
small spatial angle $\delta \theta$. The area of reflection with a
smooth boundary can be approximated as plain region. Then reflection
does not change the geometry of the bundle, and the proof for this
situation remains the same as for the bundle of trajectories with no
reflections. Hence, $p(y|x)$ exists for all $x,y \in Q$.

For convex bodies, the positivity of $p(y|x)$ is obvious, all the points
are reachable by a trajectory with no reflections.

The symmetry of the probability density function follows from the
uniformity of the distribution of the directions and reversibility
of a billiard trajectory due to the reflection law: the angle of
incidence is equal to the angle of reflection. Therefore, all the assumptions
on $p(y|x)$ are satisfied, and the distribution of points $x^i$
generated by the BW algorithm tends to the uniform distribution on $Q$.
$\Box$

\begin{theorem}
Assume $Q$ is connected, bounded and open set, the boundary of  $Q$
is piecewise smooth and for all $x,y\in Q$ there exists a
piecewise-linear path such that it connects $x$ and $y$, lies inside
$Q$ and has no more than $B$ linear parts. Then the distribution of
points $x^i$ generated by the BW algorithm tends to the uniform
distribution on $Q$ in the same sense as in Theorem 1.
\end{theorem}

\textit{Proof.} Again, the algorithm is well defined: with
probability one a point $x^{i+1} \neq x^i$ is found for arbitrary
$x^i\in Q$.

All the constraints on $Q$ are important. Connectedness guarantees
that, starting from any point, we can reach a measurable neighborhood
of any other point in $Q$. Boundedness is necessary to define the
uniform distribution on $Q$ and to prevent the trajectories to go to
infinity. Openness allows us to connect any two points with a tube
of nonzero measure. Hence, there exists a piecewise linear trajectory
connecting two arbitrary points.

Consider $p_N(y|x)$, the probability density function of $\mathbf{P}_N(A|x)$.
The inequality $p_N(y|x)>0$ holds for all integer $N>B$. The equality
$p(x^{i+1}|x^i)=p(x^i|x^{i+1})$ (reversibility)
 holds for every pair of consecutive points due to the reflection
law: the angle of incidence is equal to the angle of reflection.
Therefore, $p_N(y|x)=p_N(x|y)$.

Hence, the distribution of the subsequence $x^{0}, x^N, x^{2N},\dots$,
tends to the uniform one for $N>B$. The same is true for every subsequence
$x^{i}, x^{N+i}, x^{2N+i},\dots$. Since all the subsequences
have asymptotically uniform distribution, the distribution of points
$x^i$ generated by the BW algorithm tends to the uniform distribution on $Q$. $\Box$

% WRITE SMOOTHLY

There exist plenty of nonconvex domains that satisfy the conditions
of Theorem 2. For instance, an estimate of $B$ for the toroid is
given in Subsection 4.8. Note that the constant $B$ characterizes
the geometry of $Q$.

%In view of Theorem 1 in \cite{Smith} it suffices to prove that
%$p(y|x)>0$ for all $x,y\in\text{Int}Q$ and $p(y|x)=p(x|y)$, where
%$p(y|x)$ is probability density of transition from $x$ to $y$.

%$p(y|x)$ is not one step probability density. $p(x|y)>0$
% implies that there exists a sequence ${x^{i}}$, $i = 1, \dots K$, $K\leq B$ produced
%by BW algorithm such that $x^{1} = x$, $x^{K} = y$.

%Inequality $p(y|x)>0$ is guaranteed because all the directions are
%possible, $Q$ is connected and open, and probability of any length
%$\ell$ is positive.

\section{Discussion}

We discuss some implementation issues.

\subsection{Nonsmooth boundary points}

\begin{figure}
\begin{center}
\includegraphics[width=6.6cm]{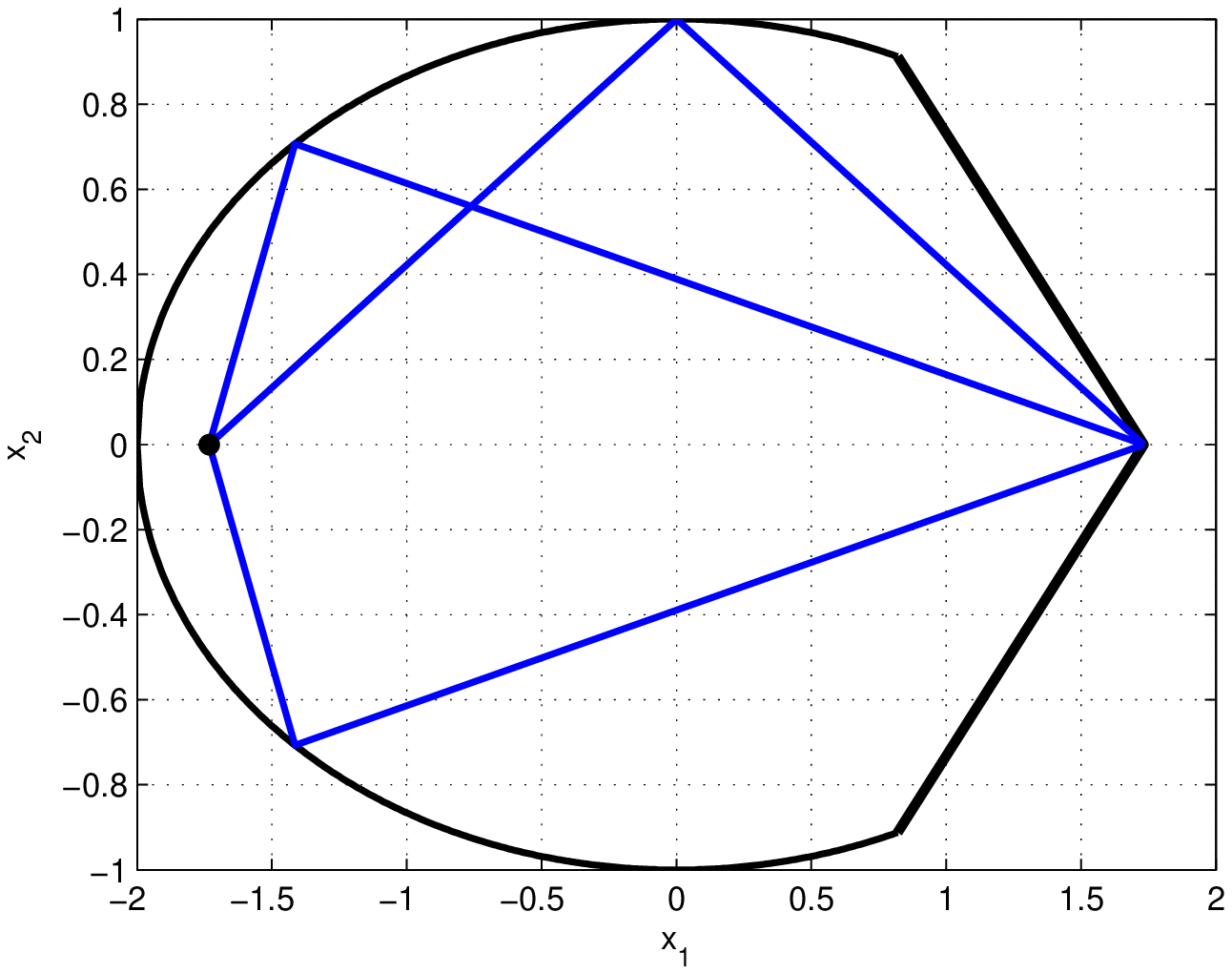}
\includegraphics[width=6.6cm]{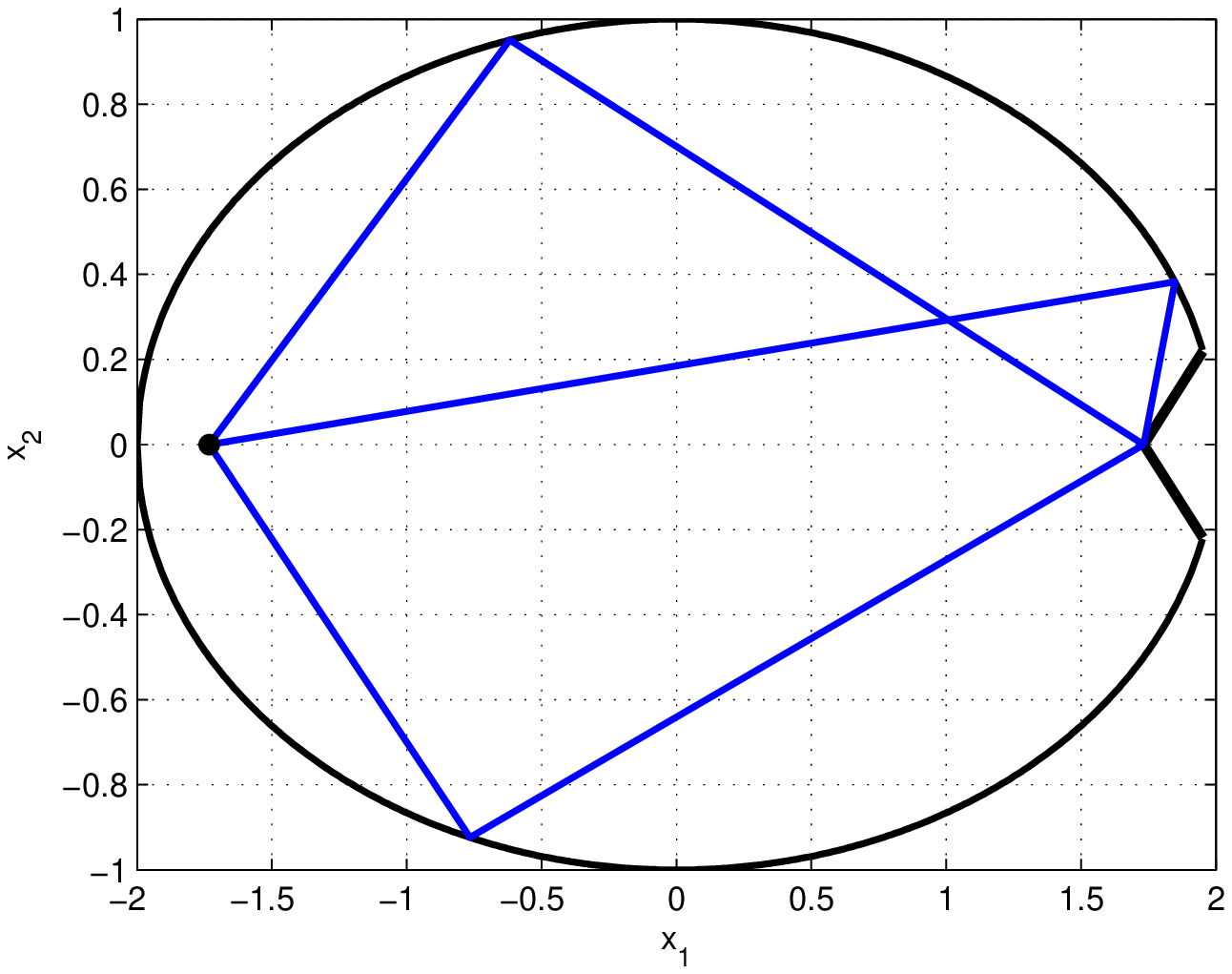}
\caption{Convex and nonconvex domains where a nonsmooth boundary point
can be achieved.}\label{fig:badshape}
\end{center}
\end{figure}

The measure of the set of points belonging to nonsmooth boundary is
zero but the probability of hitting the nonsmooth part of the
boundary is nonzero for some starting $x^0$. For instance, consider
two similar convex and nonconvex sets
\begin{gather*}
Q_1 = \{x\in \mathbb{R}^2: x_1^2/4+ x_2^2 <1, x_1 < \sqrt{3} - |x_2|\},\\
Q_2 = \{x\in \mathbb{R}^2: x_1^2/4+ x_2^2 <1, x_1 < \sqrt{3} + |x_2|\},
\end{gather*}
both being truncated ellipses (Fig. \ref{fig:badshape}). A large
portion of directions makes trajectories of length $\ell\geq 4$
starting at the focus $x^0 = (-\sqrt{3}; 0)$ hit the nonsmooth
boundary at the second focus $x^1 = (\sqrt{3}; 0)$. The reason is
that we take a particular starting point. The measure of ''bad''
starting points is zero and this effect never happens when the
starting point is taken randomly with some distribution.

\subsection{Choice of $\tau$ and $R$}

To run the algorithm we need to specify the parameters $\tau$ and $R$.
The value of $\tau$ strongly affects the behavior of the method.
For $\tau$ small enough, BW becomes slower than HR; it behaves as
a ball walk with radius $\tau$. Empirical observations show that
fast convergence to the uniform distribution is achieved for $\tau
\approx \text{diam}Q$, where $\text{diam}Q$ is the diameter of the
set $Q$.

We restrict the number of reflections by $R$ for every trajectory
(step 4 of the Algorithm). The goal is to avoid situations when
the trajectory length remains less than $\ell$ after a large number of
reflections (a typical example is addressed in Subsection 4.4). The
choice of $R$ is mostly focused on eliminating computationally hard
trajectories. The value of $R$ should be large enough to implement most of the
trajectories. But $R$ also depends on $\tau$. The longer the trajectory
one needs to implement, the more reflections are required. We usually
take $R = 10n$ to make it dimension dependent.
% Note that for certain $Q$ it may require $R>10n$ to sample efficiently.

\subsection{ Preliminary transformation of $Q$}

If $Q$ is ``ill-shaped,'' sometimes it can be improved with its
linear transformation. For instance, if $Q$ is a box $Q=\{x\in
\mathbb{R}^n: |x_i| < a_i, i = 1, \dots,n \}$, and it is far from being
cubic ($\min a_i / \max a_i \ll 1$), a simple scaling transforms $Q$
into a cube. A similar scaling transforms an ellipsoid into a ball. In the
general case ,the following scaling can be helpful. Assume $Q$ has a
barrier $F(x)$ defined on $Q$ such that $F(x) \rightarrow +\infty$
as $x \rightarrow \partial Q$. In \cite{NestNem_book94}, a special
class of self-concordant barriers is considered. For instance, for
the polytope defined by $m$ linear inequalities $Q=\{x\in
\mathbb{R}^n: (a_i,x) < b_i, i = 1,\dots, m\}$, this barrier is
$F(x)=-\sum_i \log (b_i-(a_i,x))$. Then it is easy to find
 an approximate minimum $x^*$ of $F(x)$. \emph{Dikin ellipsoid}
$E=\{x: (H(x-x^*),(x-x^*))\le 1, H= \nabla^2 F(x^*)\}$, lies in $Q$
and it is a good approximation of the polytope $Q$. Hence we can
calculate the linear mapping $T=H^{-1/2}$; by generating directions
$d'=Td$, where $d$ is uniformly distributed over the unit sphere, we
can strongly accelerate the convergence. However sometimes none of the
transformations can improve the shape of the set; the simplex is known
to be the worst-case example.

\subsection{Boundary oracle and normals}

Both the HR and BW algorithms require computation of the intersections of a
straight line (defined by the point $x^k$ and the direction $d$ of the
trajectory) with the set $Q$. We call \emph{Boundary Oracle} (BO)
the procedure that calculates the boundary of the segment
$[\underline{t}, \overline{t}]$, where
$$\underline{t}=\max\limits_{t<0} \{t: x^k+td\in \partial Q\},\quad
\overline{t}=\min\limits_{t>0} \{t: x^k+td\in \partial Q\}$$ (we
assume that $Q$ is convex, otherwise the point of the first
intersection of the straight line and the boundary of $Q$ is taken).
Thus HR needs both  $\underline{t}$ and $\overline{t}$ for every
iteration, and the computational cost of HR is equal to two BO per sample.
BW takes $\overline{t}$ for every reflection and the computation
cost of BW is one BO per reflection. In most applications, finding BO
is not a problem. For instance, if $Q$ is a polytope defined by $m$
linear inequalities
$$Q=\{x\in R^n: (a^i,x) < b_i, ~i=1,\dots, m\}$$
then BO $[\underline{t}, \overline{t}]$ can be written explicitly. Calculate
$
t_i = \frac{b_i - (a^i, x^k)}{(a^i,d)}, i = 1,\dots, m,$ and take
$$\underline{t} =\max\limits_{i:~t_i<0} t_i, \quad \overline{t} =
\min\limits_{i:~t_i>0} t_i.$$

Numerous examples of BO for other sets $Q$ (for instance, defined by
\emph{Linear Matrix Inequalities}) can be found in
\cite{PolGry:IFAC08,PolGry_AOR:09,DaShPo:SIAM10}.

Billiard walks also require the calculation of normals $s$ at the boundary
points. In most applications it is not hard; for instance, for a
polytope we have $s=a_i$, where $i$ is the index for which the maximum or
the minimum in the formulas above is achieved.

\subsection{A Comparison of HR and BW}
Our goal in the test examples below is to compare HR and BW. We use
several tools for this purpose. Sometimes theoretical considerations
can help to compare the number of iterations to quit a corner. It is
well known that HR may require too many iterations to get out of the
corner, see estimates in \cite{LovVempala_corner}. We will show that
estimates for BW are much more optimistic for many particular
examples. On the other hand, we use simulation for the comparison as
well. We exploit different tools to demonstrate that one sampling
set is closer to uniform than another. Sometimes graphical figures
in the 2D plane are quite evident. In other cases we demonstrate strong
serial correlation in the samples. Finally, we use a parametric partition
of $Q$ and compare the number of empirical frequencies with the
theoretical number for the uniform distribution via the $\chi^2$ criterion.

To make final conclusions on the comparison of the two methods, we
should have the following in mind. Of course, computationally BW is
harder than HR. It requires more BO calculations, each reflection at
the boundary also requires extra calculations for normals. We
characterize the computational complexity by the number of calls to
the BO and compare the outcomes of HR and BW obtained from the same
number of BO (the number of samples is different in this case). The
observed acceleration of convergence to the uniform distribution
often makes BW preferable to HR.

\medskip
\section{Test sets and simulation}

Some sets below are unbounded; we present them to analyze the behavior
at a corner. We say that a trajectory quits the corner if it goes to infinity.

\subsection{Plane angle}

Let the angle $Q\subset R^2$ be equal to $\alpha<\pi$. Then any
billiard trajectory quits $Q$ after no more than
$N^*=\lceil\pi/\alpha\rceil$ reflections for all initial points and
initial directions; here $\lceil a\rceil$ stands for the smallest
integer greater than or equal to~$a$. The proof of this fact is as follows: if we reflect
the angle $N$ times around its side, billiard trajectory becomes the
straight line. Every intersection of the line and the angle side
corresponds to the reflection of the billiard trajectory. A straight
line cannot intersect any straight line (not coinciding with itself)
twice, and the total number of intersections with the reflected
angle sides is no more than $N^*$. Thus $N^*$ reflections are enough
to quit the corner.

For HR we quit $Q$ with probability $1-(1-\alpha/\pi)^N$ after no
more than $N$ iterations. For $N=N^*$ being large enough, HR quits
$Q$ with probability $1-1/e=0.63$, while BW quits with probability one.

\begin{figure}
\begin{center}
\includegraphics[width=8.4cm]{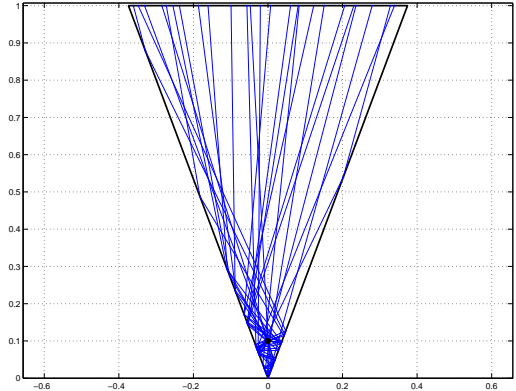}
\caption{25 trajectories reaching the line $x_2 = 1$ starting from  $[0;
0.1]$.}\label{fig:angle}
\end{center}
\end{figure}

It is of interest to estimate the average number of reflections
(over random initial directions). Consider the triangle $Q =\{x\in
\mathbb{R}^2: |x_1| \leq \text{atan} \frac{\alpha}{2}, x_2 \leq 1\}$
with one of the angles equal to $\alpha$. Let the BW trajectories start at
$x^0 =(0; 0.1)$ and calculate the number of reflections until the
trajectory reaches the line $x_2 = 1$ (i.e. quits the corner). Figure
\ref{fig:angle} depicts 25 trajectories plotted for $\alpha = \pi/4$.

The results of 5000 runs and various $\alpha$ are given in Table 1.
The empirical observations show that the average number of reflections for BW is
equal to $N^*/2$. For HR we calculate the number of iterations
until BO reaches the line $x_2 = 1$.

\begin{table}[!h]
\begin{center}
\begin{tabular}{l|c|c}
  $\alpha$ & BW & HR \\ \hline
  $\pi /2$ & 2.28 (0.87) & 2.37 (1.74) \\
  $\pi /4$ & 3.08 (1.3)  &  3.75 (2.98)\\
  $\pi /10$ & 5.94 (2.93) &  8.23 (7.1)\\
  $\pi /50$ & 25.08 (14.46) & 39.25 (34.54)
\end{tabular}
\label{tab:angle} \caption{The mean and standard deviation (in
parentheses) for the number of BW reflections and the number of HR
iterations required to quit the angle $\alpha$.}
\end{center}
\end{table}

We conclude that BW is slightly more efficient that HR.

\medskip
\subsection{Multidimensional case: polyhedral cone $Q$}

For a polyhedral cone there exists a number $M$ which does not depend on the
initial data, such that any billiard trajectory quits $Q$ after no
more than $M$ reflections (see \cite{Sinai:78}, also
\cite{Tabachnikov:95}, Theorem 7.17). However $M$ depends on the
geometry of $Q$. If $M$ is large ($M>R$) then sometimes the BW algorithm gets stuck
at $x^i$. However it can be proved that BW is well defined with probability one.

\medskip
\subsection{Orthant $Q=\{x\in R^n: x > 0\}$}

It is easy to show that a billiard trajectory quits $Q$ after no more
than $n$ reflections for an arbitrary initial point and initial
direction. Indeed, for a given initial direction $d$, every reflection makes one of the
components positive (if $d>0$
componentwise, the trajectory quits $Q$). Let $I=\{i: d_i<0\}$, then
every reflection eliminates at least one negative component of $d$,
and after no more than~$n$ reflections we have $I=\emptyset$.
\smallskip

A HR trajectory quits $Q$ with probability $2^{-(n-1)}$ after a single
iteration, thus it requires approximately $2^{n-1}$ iterations to
quit $Q$ with probability $1-1/e=0.63$. The probability to quit the
orthant after no more that $n$ iterations is $2^{-(n-2)}(1-2^{-n})$
for HR and it decreases dramatically as the dimension $n$ grows.
Hence BW is much more efficient than HR for this case. Simulations
for the cube (Subsection 4.6) confirm this statement.

\medskip

All these results show that a polyhedral corner is not a problem for
BW in contrast to HR, where the distance of the initial point to the
corner and the size of the angle plays a significant role. The results
can be extended to curvilinear corners with nondegenerate linear
approximation, i.e. if a linear approximation of a corner is a
polyhedral cone with nonempty interior.

\medskip

\subsection{Concave corner}
\begin{figure}[!h]
\centerline{\includegraphics[height=7cm]{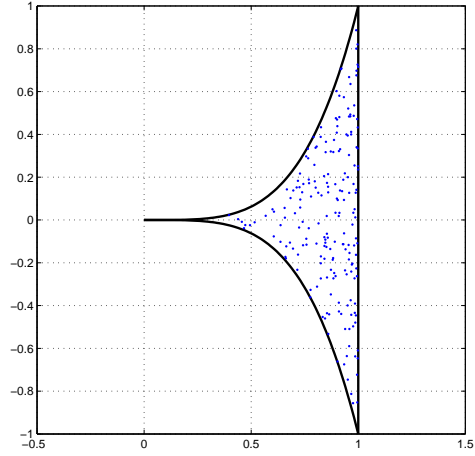}} \caption{200
points generated by BW for domain (\ref{eq:x4})} \label{fig:cusp200}
\end{figure}

In concave corners (that is, corners with concave boundaries) pure
billiard trajectories may expose a large number of reflections
\cite{Sinai:70}. Consider a typical domain (Fig. \ref{fig:cusp200})
with concave angle
\begin{equation} Q=\{x\in R^2: -x_1^4 < x_2 < x_1^4, \quad x_1\ge 1
\}.
\label{eq:x4}
\end{equation}
For a fixed $\ell$, the length of a billiard trajectory may remain
less than $\ell$ after a large number of reflections. Indeed, start the
trajectory at the point $x^0 = (0.9; \varepsilon)$, $\varepsilon$
being small enough, fix $\ell = 1$, $d = (-1; 0)$ and compute the
number of reflections required. The results are shown in Table~2. As
one can notice, the number of reflections increases dramatically as
the first coordinate of $x^0$ tends to zero, and even for $x_1^0 =
10^{-4}$, the trajectory can not be implemented. To avoid these
situations we restrict the number of reflections by $R$ in the BW
algorithm. But, in general, these ``bad'' directions are rare. Figure
\ref{fig:cusp200} depicts 200 points for domain (\ref{eq:x4}), the
average number of reflections per point is six.
\begin{table}[!h]
\begin{center}
\begin{tabular}{r|c}
  $\varepsilon$ & Number of reflections \\ \hline
    1e-3 & 746 \\
    5e-4 & 1851 \\
    4e-4 & 2480 \\
    3e-4 & 3617 \\
    2e-4 & 6158 \\
    1.1e-4 & 13496 \\
    1.01e-4 & $>5\cdot10^{6}$ \\
\end{tabular}
\label{tab:x4} \caption{The number of reflections required to implement the
trajectory of length $1$ for domain (\ref{eq:x4}) starting at $x^0 =
(0.9; \varepsilon)$ in the direction $d = (-1; 0)$.}
\end{center}
\end{table}

%\begin{figure}[!h]
%\centerline{\includegraphics[width=7cm, height=7cm]{Good.eps}
%\includegraphics[width=8.5cm, height=7cm]{Good_zoom.eps}}
%\caption{Left: the trajectory for the set (\ref{eq:cusp}) starting
%in $[-0.001; 0]$ in direction $[0; -1]$. Right: zoom of the lower
%part. }\label{fig:cusp}
%\end{figure}

%Both these situations are met at a cusp --- a curvilinear corner
%with concave boundaries with coinciding tangent planes. Typical
%example: $Q=\{x\in R^2: x_1\ge 0, ||x-a_i||\ge 1, a_1=(0,1),
%a_2=(0,-1)\}$. Then point $0$ is the attraction point for billiard
%trajectories in its neighborhood, number of reflections tends to
%$\infty$, while the length remain bounded.

\subsection{Strip}

For domains of the form $Q=\{x\in R^2: 0 < x_2 < 1, |x_1| < M\}$, $M$
being large enough, HR and BW demonstrate different abilities to
walk along $x_1$. Below we show that if one counts the average
number of steps per one BO call, BW is
approximately 6 times faster. For a random line intersecting the
strip, let $\overline{\Delta x}$ be the horizontal component of the
intersection averaged over directions. HR takes 2 BO per step, and
the average shift along $x_1$ for uniformly distributed initial
point is
$
\overline{\Delta x}_{HR} = \frac{1}{2}\overline{\Delta x}
\int\limits_0^1\int\limits_0^1 |x_1 - x_2| dx_1 dx_2  = \frac{1}{6} \overline{\Delta x}
$
per one BO.
BW gives $\overline{\Delta x}_{HR} = 1/2 \overline{\Delta x}$ for
the first reflection (1 BO) and then $\overline{\Delta x}$ for every subsequent BO.
The average shift along $x_1$ produced by BW after $N$ BO is
$
\overline{\Delta x}_{BW} = \left(1 - \frac{1}{2N}\right)\overline{\Delta x}
$
per one BO. Thus BW is 6 times more efficient than HR.

\medskip

\subsection{Cube}

For the unit cube $Q = \{x\in \mathbb{R}^n: 0< x < 1 \}$
(the inequality is understood component-wise), we can compute the next point
of the BW algorithm explicitly.

At the current point $x$, for given $\ell$ and $d$ calculate $k_i =
\lfloor x_i + \ell d_i\rfloor$ ($\lfloor x \rfloor$ is the maximal
integer less than or equal to $x$) and walk to $y$:
$$
y_i = \left\{\begin{array}{ll}
                x_i + \ell d_i - k_i, &  k_i \text{ is even}\\
                1 - (x_i + \ell d_i - k_i), &  k_i \text{ is odd}
              \end{array}
\right., \quad i = 1,\dots,n.
$$
Of course there is no need to apply MCMC algorithms for random
sampling in the cube, one can generate a vector of $n$ independent
uniform random variables over $[0, 1]$. Moreover, the shape of the
cube is so nice that the distribution of the HR points converges to the uniform
fast enough. Nevertheless it is of interest to compare BW and HR for
this simplest case.

For various dimensions $n$ we sample $N_{BW} = 1000$ points by the
BW algorithm and calculate the amount $N_{BO}$ of the BO calls needed.
Then we sample $N_{HR} = \lceil N_{BO} /2\rceil$ points by HR.
In implementing BW, we take $\tau = \sqrt{n}$, $R = 10n$, starting point
is uniform random for both sampling algorithms.

First we examine serial correlation for points produced by different
samplers. To judge about serial correlation, we partition the unit
cube into $q_i$, $i = 1,\dots, 2^n$, small cubes of equal volumes.
Then we calculate the empirical probability to proceed between
different parts $\mathbf{P}(x^{i+1} \notin q_j | x^i \in q_j)$.
Table~3 shows the results as compared to the theoretical probability $U$ for
independent uniformly distributed points (which is $2^{-n}$). One
can see that serial correlation is much stronger for HR than for the BW
samples.
\begin{table}[!h]
\begin{center}
\begin{tabular}{c|ccc}
  $n $ & BW & HR & U \\ \hline
   10  & 0.098 & 0.609 & $9\cdot 10^{-4}$ \\
   25 &  0.043 & 0.612 & $2\cdot 10^{-8}$ \\
   50  & 0.024 & 0.617 & $9\cdot 10^{-16}$
\end{tabular}
\label{tab:cube_prob} \caption{Empirical probability to proceed
between different parts of the cube for BW and HR, and the uniform
distribution.}
\end{center}
\end{table}

Then we make the $\chi^2$ frequency test for 10~000 HR points in
$\mathbb{R}^{10}$. We take 10 equal volume slabs in the $i$th
coordinate direction for $i = 1,\dots, 10$, and make 10 $\chi^2$
tests all together. The results are shown in Table 4.
\begin{table}[!h]
\begin{center}
\begin{tabular}{r|rrrrrrrrrr|r}
 & 1 & 2 & 3 & 4 & 5 & 6 & 7 & 8 & 9& 10 & $\chi^2$ st. \\ \hline
  1 &    927   &  1087 &        1096   &      985    &     987  &       992   &      963      &   979   &     1000   &      984 & 24.64\\
  2 &        1129     &    969     &    884     &   1026   &     1049   &      963     &    935  &       959& 983     &   1103 & 52.31\\
  3 &  1135 &  1134  & 970 &    951 &     983 &    976 & 980  &        961 &       822 &     1088 &  81.7 \\
  4 & 1008    &    1029     &    977       & 1046   &      961   &      870      &   932   &      971 &       1117     &   1089 & 49.05\\
  5 &      820   &     1004    &    1092 &       1043   &      956  &       960     &   1107    &    1174     &    916       &  928& 100.23\\
  6 &    1001   &     1068    &     992    &    1014   &     1051    &    1004    &     944  &       958     &    935       & 1033 & 17.71\\
  7 &     1015    &     890     &    905  &       916   &     1010     &    953  &      1028     &    982   &     1077     &   1224 & 87.93\\
  8 &      1130    &    1098      &  1074   &     1078 &       1032 &       1021    &     824   &      956    &     886    &     901 & 95.24\\
  9 &        913     &    983     &    980 &       1059    &    1023   &      902    &    1017   &     1050 &       1011    &    1062 & 28.63\\
  10 &      1056   &     1013    &     966    &     972     &    950    &    1002 &        951       &  979       & 1016      &  1095 & 19.87
  \end{tabular}
\label{tab:cube_chisquHR} \caption{The observed frequency of the HR points
in the slab $j$ in the $i$th coordinate direction and the $\chi^2$ statistics.}
\end{center}
\end{table}

We start BW with computational complexity 20 000 BO and obtain 2148
points ($\sim215$ per slab). The $\chi^2$ frequency test results are
shown in Table 5.
\begin{table}[!h]
\begin{center}
\begin{tabular}{r|rrrrrrrrrr|r}
 & 1 & 2 & 3 & 4 & 5 & 6 & 7 & 8 & 9& 10 & $\chi^2$ st. \\ \hline
1 & 198  & 204  & 188 &  218 &  216  & 224  & 218 &  242 &  191 &  249 & 17.22 \\
2 &   234  & 207 &  196 &  226 &  204  & 196  & 225 &  230  & 220 &  210 & 8.21\\
3 &    210  & 230 &  218 &  201 &  202 &  214 &  232 &  200&   222&   219 & 5.6\\
4 &   242  & 243 &  203 &  198 &  202  & 232 &  221  & 208 &  220 &  179 & 17.7\\
5 &   211  & 231 &  184 &  236 &  229  & 206 &  210  & 235 &  192 &  214 & 13.52\\
6 &   209  & 193 &  242 &  205 &  216  & 208 &  212  & 223 &  222 &  218 & 7.12\\
7 &   190  & 223 &  226 &  233 &  197  & 217 &  226  & 195 &  200 &  241 & 13.42\\
8 &   200  & 231 &  199 &  191 &  207  & 211 &  212  & 220 &  247 &  230 & 12.46\\
9 &   231  & 213  & 212 &  224 &  189  & 234 &  209  & 225 &  197 &  214 & 8.6\\
10 &   204 &  237 &  227 &  198 &  201 &  230 &  208 &  215 &  211 &  217 & 7.11\\
\end{tabular}
\label{tab:cube_chisquBW} \caption{The observed frequency of the BW
points in the slab $j$ in the $i$th coordinate direction and the
$\chi^2$ statistics.}
\end{center}
\end{table}

Upper and lower $\chi^2$ values for $10\%$ statistical significance for 9 degrees of freedom are $[3.3, 16.9]$ (for two tailed $\chi^2$ test).
Thus HR fails all 10 $\chi^2$ tests while BW fails just 2 out of the 10 tests.

\subsection{Simplex}

The next test set is the standard $n$-dimensional simplex
$$
Q = \{ x_i > 0, \sum x_i = 1, i = 0, 1, \dots, n\}.
$$

The simplex is a set with many corners and the geometry of simplex cannot be improved by any affine transformation.
We know that for HR
walk it takes a lot of iterations to get out of a corner, thus it is
interesting to compare HR and BW.

Smooth boundary of $Q$ is specified by the points $\partial Q =\{x\in
\mathbb{R}^{n+1}: x_k=0, x_i \neq 0, i = 0, \dots, n, i \neq k\}$,
and the internal normal of unit length for these points is
$$
s =\sqrt{\frac{1}{n(n+1)}}\left[-1, \dots, \underbrace{n}_k,
\dots,-1\right]^T.
$$
The length of any edge of the simplex is $\sqrt{2}$ for every
dimension $n$, the choice of the parameter $\tau = \sqrt{2}$ seems
reasonable.

%Fig. \ref{fig:tri} shows 300 points generated by HR
%(black) and BW (blue) for standard 2-simplex.
%\begin{figure}[!h]
%\centerline{\includegraphics[width=11cm,
%height=10cm]{triangle_R3_2.eps}} \caption{300 points generated by HR
%(black) and BW (blue) points for standard 2-simplex.}\label{fig:tri}
%\end{figure}

\begin{figure}[!h]
%\centerline{\includegraphics[width=9cm,
%height=7.5cm]{S2_30_many10.eps}
%\includegraphics[width=9cm, height=7.5cm]{S2_300_many10.eps}}
\centerline{\includegraphics[width=9.5cm]{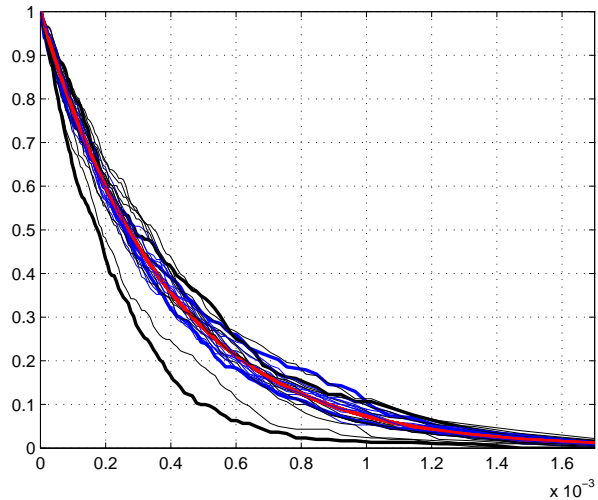}}
\caption{The portion of points contained in $S_{\alpha}$ for uniformly
distributed points (red), HR (black) and BW (blue). $n=50$, 300
points. The horizontal axis corresponds to the parameter
$\alpha$.}\label{fig:cdf}
\end{figure}

Note that for $n=2$, the samples visually look uniformly distributed for
both algorithms. To decide about uniformity more rigorously in the
multidimensional case, consider the sequence of enclosed simplices
$S_{\alpha} = \{x\in \mathbb{R}^{n+1}: x_i \geq \alpha, \sum x_i =
1\}$, $0 \leq \alpha\leq \frac{1}{n+1}$. For $\alpha = 0$, $S_0$ is
the initial simplex, and for $\alpha = \frac{1}{n+1}$, the simplex
$S_{\alpha}$ contains one point. Let $\widehat{f}(\alpha)$ be the
portion of points contained in $S_{\alpha}$, and denote $ f(\alpha)
= \text{vol}S_{\alpha}/\text{vol}S_0 = (1 - (n+1)\alpha)^n. $ Figure
\ref{fig:cdf} shows $\widehat{f}(\alpha)$ for $n= 50, N=300, x^0 =
\{1/(n+1), \dots 1/(n+1)\}$. The red line corresponds to the uniformly
distributed points, the black line describes the distribution for the HR
points, and the blue line for BW points. We conclude that for BW samples, the
empirical values of $\widehat{f}(\alpha)$ are much closer to the mean
value $f(\alpha)$ than for the HR samples.

%Another advantage of BW --- its ability to quit a corner fast
%enough. We take a specific starting point
% close to the corner $x^0 = \left[1 - \varepsilon,
%\frac{\varepsilon}{n}, \dots, \frac{\varepsilon}{n}\right]^T,
%\varepsilon = 0.1$.

%Fig. \ref{fig:corner_norm} depicts the distance between $x^i$ and
%the vertex of the simplex $c = [1, 0, \dots, 0]^T$. HR fails  to get
%out of the corner, 100 samples remain in the neighborhood of the
%starting point. BW demonstrate strong mixing, after a few iterations
%points look like uniformly distributed ones.

%\begin{figure}[!h]
%\centerline{\includegraphics[width=9cm, height=7cm]{S50_200bad.eps}
%} \caption{Portion of points contained in $S_{\alpha}$ for uniformly
%distributed points (red), HR (black) and BW (blue). $n=50$, $N=200$,
%$x^0 = \left[1 - \varepsilon, \frac{\varepsilon}{n}, \dots,
%\frac{\varepsilon}{n}\right]^T$.}\label{fig:cdf_bad}
%\end{figure}

%\begin{figure}[!h]
%\centerline{\includegraphics[width=9.5cm]{corner50_1.eps}}
%\caption{Distance from the corner $[1, 0, \dots, 0]^T$ for the first
%100 points of walks in 50-simplex for $x^0 = \left[1 - \varepsilon,
%\frac{\varepsilon}{n}, \dots, \frac{\varepsilon}{n}\right]^T$,
%$\varepsilon = 0.1$. Uniformly distributed points (red), HR (black)
%and BW (blue).}\label{fig:corner_norm}
%\end{figure}

We also perform two $\chi^2$ tests for $n = 10$. For the first one
we partition $Q$ into 10 simplices $Q = S_{0} \supset S_{\alpha_1}
\supset \dots \supset S_{\alpha_{10}} = \emptyset$ such that
$\text{Vol}(S_{\alpha_i} \setminus S_{\alpha_{i+1}}) = \text{Vol}
S_{\alpha_9} = \frac{1}{10} \text{Vol}S_{0}$. These
differences $S_i \setminus S_{\alpha_{i+1}}$ are of various geometry
but their volumes are equal. For the second test we take $n+1$
subsets of the same volume and geometry, these subsets $Q_i$ contain
points mostly close to the selected vertex $v^i$: $$Q_i = \{x\in Q:
||x-v^i||_2<||x-v^j||_2, \quad j\neq i\}.$$ Restricting ourselves to 20 000 BO, we obtain
10 000 HR points and about 2 000 BW points. Tables 6 and 7 present the experimental results.
\begin{table}[!h]
\begin{center}
\begin{tabular}{r|rrrrrrrrrr|r}
 & 1 & 2 & 3 & 4 & 5 & 6 & 7 & 8 & 9& 10 & $\chi^2$  \\ \hline
HR & 855  & 917  & 925 &  897 &  997  & 978  & 1024 &  1025 &  1080 &  1303 & 144.04 \\
BW &   152  & 175 &  163 &  177 &  189  & 192  & 206 &  182  & 214 &  241 & 6.06
\end{tabular}
\label{tab:cube_chisquBW} \caption{The observed frequencies for the
HR and BW points in the subsets $S_{\alpha_i} \setminus
S_{\alpha_{i+1}}$ and the $\chi^2$ statistics.}
\end{center}
\end{table}
\begin{table}[!h]
\begin{center}
\begin{tabular}{r|rrrrrrrrrrr|r}
 & 0 & 1 & 2 & 3 & 4 & 5 & 6 & 7 & 8 & 9& 10 & $\chi^2$ \\ \hline
HR & 976 & 751 & 1050 & 1018 & 826 & 676 & 1084 & 521 & 1028 & 1424 & 947 & 697.54 \\
BW &  172 & 188 & 167 & 177 & 158 & 179 & 176 & 154 & 147 & 178 & 195 & 12.1
\end{tabular}
\label{tab:cube_chisquBW} \caption{The observed frequencies for the HR and
BW points in the subsets $Q_i$, $i = 0, \dots, 10$ and the $\chi^2$
statistics.}
\end{center}
\end{table}

Recalling the upper and lower $\chi^2$ test values $[3.3, 16.9]$,
we conclude that the BW points fit uniform distribution while the HR points
do not.

\subsection{Toroid}

Both the HR and BW algorithms are applicable to nonconvex sets.
Consider the toroid formed by an $n$-dimensional ball of radius $r$
with its center rotating over a circle in the $(x_1, x_2)$-plane:

\begin{equation}
\label{eq:torus} Q =\{x\in \mathbb{R}^n: ||x - c_x|| < r\},
\end{equation}
where ${c_x}_i = \frac{x_i}{\sqrt{x_1^2 + x_2^2}}$, $i = 1,2$,
${c_x}_i = 0$, $i>2$.

The conditions of Theorem 2 are satisfied with $B =
\left\lceil\frac{\pi}{2\arccos\frac{1-r}{1+r}}\right\rceil + 1$,
i.e. for all $x,y\in Q$ there exists a piecewise-linear path such
that it connects $x$ and $y$, lies inside $Q$, and has no more than
$B$ linear parts.

\begin{figure}[!h]
\centerline{\includegraphics[width=10.5cm]{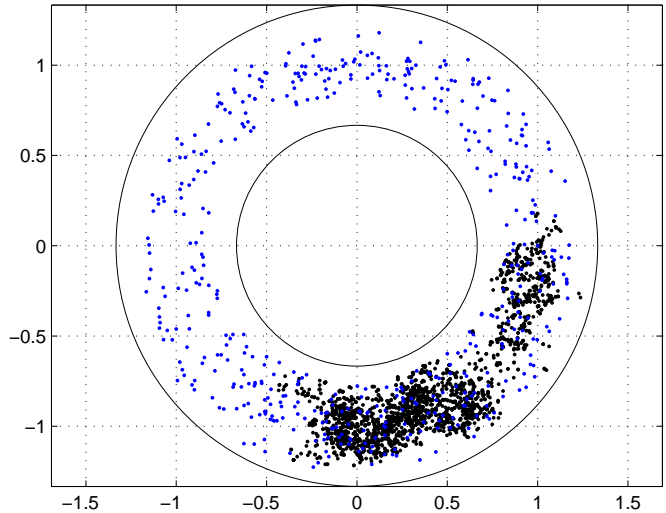}}
%\includegraphics[width=9cm, height=7cm]{Tor50_10000.eps}
%\centerline{\includegraphics[width=9cm,
%height=7.5cm]{Tor10_1000hist.eps}
%\includegraphics[width=9cm, height=7.5cm]{Tor50_10000hist.eps}}
\caption{The $(x_1, x_2)$-projection of HR points (black) and BW points
(blue) for toroid (\ref{eq:torus}). $n=10$, $N_{BW}=500$, $N_{HR} =
1764$. }\label{fig:torus}
\end{figure}

Figure \ref{fig:torus} depicts $500$ BW samples and $1~764$ HR
points (projected onto the $(x_1, x_2)$-plane) for the set
(\ref{eq:torus})  of
dimension $10$ and $r = 1/3$. The number of samples
 is different because implementation of 500 BW samples requires 1~764 BO calculations. %(left)
%and $50$ (right), samples are projected to
% $(x_1, x_2)$-plane.
HR points are plotted with black dots, BW points
with blue ones.
%The lower part of Fig.
%\ref{fig:torus}  contains histograms of the angle distribution.

It can be easily seen that the angles of the BW points are much more
uniformly distributed than those for the HR points, the latter remain in the
neighborhood of the initial point. Note that the visual lack of
uniformity in the radial direction is an ''optical effect'' because we
provide a 2D projection of the 10D picture.

%Similar behavior of random walks is observed for the set
%\begin{equation}
%\label{eq:korka} Q =\{x\in \mathbb{R}^n:\theta \leq ||x|| \leq 1,
%0<\theta <1\}
%\end{equation}

%Compare to (\ref{eq:torus}) HR demonstrates faster convergence. The
%reason is that the portion of directions that allow HR to make large
%step is greater for (\ref{eq:korka}) than for (\ref{eq:torus}). Fig.
%\ref{eq:korka} shows histograms of the angle distribution ...

\section{Applications}

In this paper we do not address numerous applications of the new version
of random sampling. We can mention just few of them: global
optimization (in particular, concave programming), control problems,
robustness issues, numerical integration, calculation of the volume and
the center of gravity and so on; see, for instance, our previous
papers
\cite{PolGry:IFAC08,PolGry_AOR:09,DaShPo:SIAM10,PolGry:MSC10}. We
plan to consider these applications in future works.

\section{Acknowledgements}
The initial impulse for the research was given by a student
Alexander Rodin, who proposed to use physical model of gas diffusion
for random sampling. The suggestions of Yakov Sinai on billiard theory
were extremely helpful.

\bibliographystyle{elsarticle-num}
\bibliography{polgryaz}

\begin{thebibliography}{10}
\expandafter\ifx\csname url\endcsname\relax
  \def\url#1{\texttt{#1}}\fi
\expandafter\ifx\csname urlprefix\endcsname\relax\def\urlprefix{URL }\fi
\expandafter\ifx\csname href\endcsname\relax
  \def\href#1#2{#2} \def\path#1{#1}\fi

\bibitem{TeCaDa:04}
R.~Tempo, G.~Calafiore, F.~Dabbene, Randomized Algorithms for Analysis and
  Control of Uncertain Systems, Communications and Control Engineering Series,
  Springer-Verlag, London, 2004.

\bibitem{Rubin_MCbook:2008}
R.~Rubinstein, D.~Kroese, Simulation and the Monte Carlo Method, Wiley, NJ,
  2008.

\bibitem{Gilks_book96}
W.~Gilks, S.~Richardson, D.~Spiegelhalter, {M}arkov {C}hain {M}onte {C}arlo,
  Chapmen and Hall, 1996.

\bibitem{Diaconis09}
P.~Diaconis, The markov chain monte carlo revolution, Bull. of the AMS 46~(2)
  (2009) 175--205.

\bibitem{DyerFriKannan:91}
M.~Dyer, A.~Frieze, R.~Kannan, A random polynomial-time algorithm for
  approximating the volume of convex bodies, Journal of the ACM 38~(1) (1991)
  1--17.

\bibitem{LovaszSomonovits:93}
L.~Lovasz, M.~Somonovits, Random walks in a convex body and an improved volume
  algorithm, Random Structures \& Algorithms 4~(4) (1993) 359--412.

\bibitem{LovaszDeak:12}
L.~Lovasz, I.~Deak, Computational results of an ${O}^{*}(n^4)$ volume
  algorithm, European Journal of Operational Research 216 (2012) 152--161.

\bibitem{Turchin_eng}
V.~Turchin, On the computation of multidimensional integrals by the {M}onte
  {C}arlo method, Theory of Probability and its Applications 16~(4) (1971)
  720--724.

\bibitem{Smith}
R.~Smith, Efficient {M}onte {C}arlo procedures for generating points uniformly
  distributed over bounded regions, Operations Research 32~(6) (1984)
  1296--1308.

\bibitem{PolGry:IFAC08}
B.~Polyak, E.~Gryazina, Hit-and-{R}un: New design technique for stabilization,
  robustness and optimization of linear systems, in: Proceedings of the IFAC
  World Congress, Seoul, South Korea, 2008, pp. 376--380.

\bibitem{PolGry_AOR:09}
B.~Polyak, E.~Gryazina, Randomized methods based on new {M}onte {C}arlo schemes
  for control and optimization, Annals of Operational Research 189~(1) (2011)
  343--356.

\bibitem{DaShPo:SIAM10}
F.~Dabbene, P.~Shcherbakov, B.~Polyak, A randomized cutting plane method with
  probabilistic geometric convergence, SIAM Journal of Optimization 20~(6)
  (2010) 3185--3207.

\bibitem{Tervonen_etal:EJOR:2013}
T.~Tervonen, G.~van Valkenhoef, N.~Basturk, D.~Postmus, Hit-and-run enables
  efficient weight generation for simulation-based multiple criteria decision
  analysis, European Journal of Operational Research 224 (2013) 552--559.

\bibitem{NestNem_book94}
Y.~Nesterov, A.~Nemirovsky, Interior Point Polynomial Methods in Convex
  Programming, SIAM, Philadelphia, 1994.

\bibitem{PolGry:MSC10}
B.~Polyak, E.~Gryazina, {M}arkov {C}hain {M}onte {C}arlo method exploiting
  barrier functions with applications to control and optimization, in: IEEE
  Multi-Conference on Systems and Control, 2010, pp. 1553--1557.

\bibitem{Tabachnikov:95}
S.~Tabachnikov, Geometry and Billiards, RI: Amer. Math. Soc., 1995.

\bibitem{GalpZem:90_billiards}
G.~Galperin, A.~Zemlyakov, Mathematical Billiards, Nauka, Moscow (in Russian),
  1990.

\bibitem{Sinai:70}
Y.~G. Sinai, Dynamical systems with elastic reflections, Russian Mathematical
  Surveys 25~(2) (1970) 137--189.

\bibitem{Sinai:78}
Y.~G. Sinai, Billiard trajectories in a polyhedral angle, Russian Mathematical
  Surveys 33~(1) (1978) 219--220.

\bibitem{Kozlov:billiards_book:91}
V.~V. Kozlov, D.~V. Treshchev, Billiards: A Genetic Introduction to the
  Dynamics of Systems with Impacts, Vol.~89, Translations of Mathematical
  Monographs by American Mathematical Society, Providence, RI, 1991.

\bibitem{evans2001stochastic}
S.~N. Evans, Stochastic billiards on general tables, Annals of Applied
  Probability (2001) 419--437.

\bibitem{Shake-n-Bake:91}
C.~Boender, R.~Caron, J.~McDonald, A.~R. Kan, H.~Romeijn, R.~Smith, J.~Telgen,
  A.~Vorst, Shake-and-bake algorithms for generating uniform points on the
  boundary of bounded polyhedra, Operations research 39~(6) (1991) 945--954.

\bibitem{dieker2013StochBill}
A.~Dieker, S.~S. Vempala,
  \href{http://www2.isye.gatech.edu/$\sim$adieker3/publications/
  stochasticbilliard-submit.pdf}{Stochastic billiards for sampling from the
  boundary of a convex set}.
\newline\urlprefix\url{http://www2.isye.gatech.edu/$\sim$adieker3/publications/
  stochasticbilliard-submit.pdf}

\bibitem{LovVempala_corner}
L.~Lovasz, S.~Vempala, Hit-and-run from a corner, in: Proceedings of the 36th
  annual ACM symposium on Theory of computing, Chicago, IL, USA, 2004, pp.
  310--314.

\end{thebibliography}

% http://mathworld.wolfram.com/Billiards.html

%\begin{thebibliography}{10}

%\bibitem{Fam} Fam A. and Meditch J. A canonical parameter space for
%linear systems design. IEEE TAC, 1978, 23(3).

%\bibitem{Jury}
%Jury E. Inners and Stability of Dynamic Systems. New York: Wiley,
%1974.

%\end{thebibliography}

\end{document}